\documentclass{llncs}

\usepackage{epsfig,epstopdf,graphicx, amsmath, amssymb, booktabs, algorithm, algpseudocode, tikz, cite, lscape, booktabs, geometry}

\begin{document}

\title{A Geometrical Branch-and-Price (GEOM-BP) Algorithm for Big Bin Packing Problems}

\author{Masoud Ataei \thanks{mataei@mathstat.yorku.ca} \and
Shengyuan Chen{}
}

\institute{Department of Mathematics and Statistics, York University, Ontario, Canada
}

\maketitle              

\begin{abstract}
Bin packing problem examines the minimum number of identical bins needed
to pack a set of items of various weights. This problem arises in various areas of the artificial intelligence demanding derivation of the exact solutions in the shortest amount of time. Employing branch-and-bound and
column generation techniques to derive the exact solutions to this problem, usually requires designation of problem-specific branching rules compatible with the nature of
the polluted pricing sub-problem of column generation. In this work, we present
a new approach to deal with the forbidden bins which handles
two-dimensional knapsack problems. Furthermore, a set of diving criteria
are introduced which emphasize the importance of the geometrical features of the bins.
It is further shown that efficiency of the column generation technique
could significantly get improved using an implicit sectional
pricing scheme. The proposed algorithm outperforms the current state-of-the-art algorithms in number of
the benchmark instances solved in less than one
minute.
\end{abstract}
\section{Introduction}
\label{Intro}
In a bin packing problem (BPP), we are given an unlimited number of
bins of identical capacity $c > 0$ and a set of items $I = \{i_1 , \dots , i_n \}$ each of which having a weight $w_i$ $(0 < w_i \leq c)$, and the goal is to pack all the items in a minimum number of bins without exceeding bins’ capacities. For
practical problems, we can assume that the weights and capacities are
integers.

The BPP arises in various areas of the artificial intelligence including \textit{cloud resource allocation and management} 
\cite{ren2017competitiveness, hallawi2017multi, liao2017resource, kaseb2017cloud, wang2017joint}, 
\textit{automatic power plant assignments}
\cite{benazouz2017safety, yao2017outage} and
\textit{multimedia processing} 
\cite{yao2017outage,gulati2017multimedia}.
In many of these applications, the exact solutions to the problem are often required to be obtained. 

To the best of our knowledge, the most recent review on exact algorithms used for solving BPP is presented by Delorme et al. 
\cite{delorme2016bin,delorme2017bpplib} 
where the authors have reviewed various formulations of BPP and conducted extensive experiments on benchmark instances. Broadly speaking, the develped exact algorithms fall into one of the following categories: branch-and-bound
\cite{martello1990lower, scholl1997bison, lysgaard2004new}, 
branch-and-price 
\cite{belov2006branch, vance1994solving, vance1998branch, ryan1981integer, vanderbeck1999computational, vanderbeck2000dantzig} 
and pseudo-polynomial algorithms 
\cite{rao1976cutting, dyckhoff1981new, stadtler1988comparison, de1999exact, cambazard2010propagating, brandao2016bin}.

In this paper, we present an exact algorithm for solving the BPP and that
relies on branch-and-bound and column generation techniques. Employing
column generation in a branch-and-bound tree (also known as branch-and-price algorithm), requires the development of rigorous methods to solve the
pricing sub-problem of the column generation. Specifically, when branching
occurs on variables of the Gilmore-Gomory model 
\cite{gilmore1961linear, gilmore1963linear}, 
bins known as
the forbidden bins (patterns) are dealt with while proceeding into the depth
of the branch-and-bound tree. The forbidden bins should then be systematically excluded from the search domain of the one-dimensional knapsack
problem (1D-KSP).

To address the issue of having a polluted pricing sub-problem in column
generation, two main approaches have been introduced in the literature.
The first approach determines the $K$-best solutions of the 1D-KSP at the
$(k-1)$-th level of the tree. This guarantees that enough number of bins is
available to be added to the restricted master problem in case all bins appear
to be forbidden ones during column generation. Although some attempts
(see \cite{leao2014determining, sarin2005schedule}) have been made to determine the $K$-best solutions of the 1D-
KSP, there still persists a lack of methods featuring universal applicability.
In the second approach, a set of branching rules is designed that avoids
emergence of forbidden bins in the pricing sub-problem of column generation
\cite{vance1994solving, vance1998branch, ryan1981integer, vanderbeck2000dantzig}. 
The drawbacks of such branching schemes have been discussed
in details by Vanderbeck
\cite{vanderbeck2011branching}.

As an alternative, a new approach to deal with forbidden bins is proposed in this work. A constraint referred to as a decrement constraint is added to
the 1D-KSP whenever one of the forbidden bins is met during column generation. This extra constraint compels the pricing sub-problem to generate the
next feasible solution of the 1D-KSP. Subsequently, the generated solution
(bin) is passed to the restricted master problem. Premature termination
of column generation might occur when employing our proposed method.
The consequences of such an undesired termination on the branch-and-price
algorithm are investigated and resolved in this paper.

Another important aspect of efficient branch-and-price algorithms has to do with primal heuristics. A thorough study of different primal heuristics
like diving, relaxation induced neighborhood search and local branching has
been reported earlier 
\cite{sadykov2015primal}. 
Yet, the importance of the geometrical features
of bins (such as Pythagorean means of weights of the items packed into bins) in diving methods has not received any attention in the literature.
In this work, we further exploit the geometrical features of bins in enhancing performance of the diving methods and propose an effective primal
heuristic that is called batch diving. Our concern, like in other primal
heuristics, is to take full advantage of LP relaxation solutions for constructing the integer solutions. Batch diving is in fact a multi-dimensional knapsack problem that plays a key role in accelerating the proposed geometrical branch-and-price (\texttt{GEOM-BP}) algorithm.

In Sect. 2, the set partitioning formulation of the BPP, column generation and implicit sectional pricing scheme are presented. In this section, a
simple method named subset-sum-$\tilde{n}$ heuristic used to initialize column generation is also proposed. Geometrical interpretations of packing problems
and their use in enhancing the diving methods is addressed in Sect. 3.
The mathematical model that considers batch diving as a generalization of
geometrical diving methods is described in Sect. 4. In Sect. 5, we
present a branching scheme that divides the search region of a problem in
a way that exploration of the nodes with forbidden bins becomes an indispensable aspect of the Geom-BP algorithm. Finally, computational results
and discussions are presented in Sect. 7.

\section{Mathematical Formulation}
\label{MathematicalFormulation}
A BPP as defined earlier could be considered as a special case of one-dimensional cutting stock problem (1D-CSP). In 1D-CSP, each item is further associated with a demand $d_i$ $(d_i \in Z^+)$. The set of items $I = \{i_1 , \dots , i_n \}$
will then contain items with unique weights only, and the items accommodated into bins of the solutions should further satisfy their corresponding
demands.

It is noteworthy that in the context of the CSP, bins are usually referred to as cutting rolls, cutting patterns or simply patterns with the latter
appellation being the most frequent in the literature. Also, when column
generation technique is involved, the terms column, bin and pattern are used
synonymously.

Even though modeling and solving either of the BPP or 1D-CSP yields
optimal number of bins, 1D-CSP formulation is preferred when iterative
matrix-based algorithms like column generation are used in the solving procedure. This is due to the dimensions of the matrices that are lower when
the cutting stock approach is adopted, considering the fact that the majority
of BPP instances from industry as well as the ones found in the literature
involve items with same weights. For the purpose of this work, the classical
definition of BPP given in Sect. \ref{Intro} is altered to include demands as extra
attributes of items in addition to their weights.

\subsection{Set Partitioning Formulation and Column Generation}
Set partitioning formulation for BPP is introduced as
\begin{subequations} 
\label{master}
\begin{alignat}{3}
\label{master_1}
\min & \ \sum_{k=1}^{K} \lambda_k  \\
\label{master_2}
\mathrm{s.t.} & \ \sum_{k=1}^{K} x_{ik} \lambda_k = d_i && \quad i=1,\dots,n \\
 & \ \lambda_k \geq{0} \ \mathrm{ and \ integer} &&\quad k=1,\dots,K \ ,
\end{alignat}
\end{subequations}
where $\lambda_k$ is the load of bin $k$ (number of times bin $k$ is used), and $K$ denotes
an upper bound on the number of possible bins. The variable $x_{ik}$ in (\ref{master_2})
is the number of replications of item $i$ in bin $k$, and this set of constraints
ensures that all the items satisfy their demands. Also, feasibility of the bins
are ensured by the constraint
\begin{equation}
\left\{ \sum_{i=1}^{n} w_i x_{ik} \leq c \ , \ k=1,\dots,K \ , \ x_{ik}\in \{0,\dots,d_i \} \ , \ i=1,\dots,n  \right\} \ .
\end{equation}

To solve the LP relaxation of model (\ref{master}), it is impractical to enumerate
all the possible bins, a fact that calls for column generation technique. By
dropping the integrality constraint in (\ref{master_2}), the master problem for the BPP
is derived as
\begin{subequations} 
\label{master_LP}
\begin{alignat}{3}
\label{master_LP_1}
\min & \ \sum_{k=1}^{K} \lambda_k  \\
\label{master_LP_2}
\mathrm{s.t.} & \ \sum_{k=1}^{K} x_{ik} \lambda_k = d_i && \quad i=1,\dots,n \\
 & \ \lambda_k \geq{0} &&\quad k=1,\dots,K \ .
\end{alignat}
\end{subequations}

Solving model (\ref{master_LP}) using column generation starts by first defining the
restricted master problem (RMP) initialized by a basic feasible solution. The RMP is then solved to optimality by the revised simplex method
and its dual prices denoted by $\boldsymbol{\pi}$, are assumed to be the cost coefficients of
the pricing sub-problem
\begin{equation}
\label{1D_KSP_bou}
\nu_{\mathrm{exact}}(\boldsymbol{\pi}) = \max \left\{\sum_{i=1}^{n} \pi_i x_{i} - 1: \sum_{i=1}^{n} w_i x_{i} \leq c \ , \ x_{i}\in \{0,\dots,d_i \} \ , \ i=1,\dots,n  \right\} \ ,
\end{equation}
and the generated bin by solving the 1D-KSP (\ref{1D_KSP_bou}) is subsequently added to
the RMP. A given bin terminates the column generation process in case its
reduced cost $\nu_{\mathrm{exact}}(\boldsymbol{\pi}) \leq 0$. Were it to be the case, termination of column
generation indicates that all of the bins with positive reduced costs have been priced
out.

It is pointed out that the method used for solving the pricing sub-problem
(\ref{1D_KSP_bou}) should be capable of finding the exact solution of the 1D-KSP problem
upon termination of column generation. Otherwise, column generation could
be terminated prematurely. In this case, the objective value of the master
problem solved at the root node of the tree would not represent a lower
bound on the BPP.

\subsection{Implicit Sectional Pricing Scheme}
There are two major factors contributing to the amount of computational
time when master problem is solved by column generation: i) the dimensions
of the basis of the master problem, ii) the dimensions of the pricing sub-problem. The dimensions of the basis determine the performance of basis-
updating techniques. This is due to the fact that most of the computational
time needed for column generation technique elapses while updating the
basis. On the other hand, dimensions of the pricing sub-problem determine
how fast the new columns (bins) are generated during column generation.

Dimensions of the pricing sub-problem could be reduced through solving
a binary 1D-KSP
\begin{equation}
\label{1D_KSP_binary}
\nu_{\mathrm{sec}}(\boldsymbol{\pi}) = \max \left\{\sum_{i=1}^{n} \pi_i x_{i} - 1: \sum_{i=1}^{n} w_i x_{i} \leq c \ , \ x_{i}\in \{0,1 \} \ , \ i=1,\dots,n  \right\} \ ,
\end{equation}
where the generated columns contain unique items (no items with same
weights). This set of columns form a section implicitly and our computational experiments show that picking the entering column from the defined
implicit section (\ref{1D_KSP_binary}) is efficient. However, there are other columns not defined
by (\ref{1D_KSP_binary}). For this reason, whenever a column generated by solving model (\ref{1D_KSP_binary})
is priced out with negative reduced cost, the original pricing sub-problem
(\ref{1D_KSP_bou}) is called to prove optimality of the column generation.

\subsection{Subset-sum-$\tilde{n}$ Heuristic}
Quiroz-Castellanos et al. \cite{quiroz2015grouping} propose a first fit algorithm with $\tilde{n}$ pre-allocated-items (FF-$\tilde{n}$) to be used as an upper bounding technique for BPP.
In FF-$\tilde{n}$, the items $\tilde{N} = \left\{i\in I : w_i \geq \lceil{\frac{c}{2}} \rceil \right\}$ are placed into separate bins,
then the first fit algorithm is ran to place the remaining items into the bins.

In our proposed subset-sum-$\tilde{n}$ heuristic, the $\tilde{n}$ large items $(\tilde{n} = |\tilde{N}|)$, similarly to FF-$\tilde{n}$, are accommodated into separate bins. Subsequently, instead
of using the first fit algorithm, the 1D-KSP
\begin{equation}
\label{subset_sum}
\max \left\{\sum_{i=1}^{n} w_i x_{i} : \sum_{i=1}^{n} w_i x_{i} \leq c \ , \ x_{i}\in \{0,\dots,d_i \} \ , \ i=1,\dots,n  \right\} \ ,
\end{equation}
also known as the subset-sum problem is employed and solved iteratively
until all the remaining items are assigned to bins.
The objective of model (\ref{subset_sum}) is to maximize the fill of the bins. Such well-filled bins could effectively be used as good candidates to initialize column
generation. 

\section{Geometrical Diving}
\label{DivingSec}
In diving methods, branching usually occurs on a bin with the highest
value in the solution of master problem and is followed by re-optimization of
the residual problem 
\cite{vanderbeck2011branching, atamturk2005integer}. This process is continued at each node until
an integer solution referred to as primal heuristic solution is achieved. 

In this work, the diving criteria are formulated using the Lehmer mean of the weights of items packed into
the bins. Recall that, for a
given set of weights $W = \{w_j: w_j \in \mathbb{R} \ , \ j = 1, 2, \dots , J\} $, the Lehmer mean function $L_p(W)$ is defined as
\begin{equation}
L_p(W) = \dfrac{\sum\limits_{j=1}^{J} w_j^{p} }{ \sum\limits_{j=1}^{J} w_j^{p-1} } \ ,
\end{equation}
for every $p\in \mathbb{R}$ \cite{bullen2013handbook}. Special Lehmer means are obtained when values of $p$ are set to $0$ and
$2$, where the terms $L_0(W)$ and $L_2(W$) are known as the harmonic and contra-harmonic means of the set $W$, respectively.

Obtaining a solution to BPP by packing the larger
items first, could then be modeled by considering the contra-harmonic mean of
the bins as diving criterion. This is due to the fact that small items make small contributions to the
contra-harmonic mean of their containing bins. Among all the bins obtained
by solving the master problem at each node, the bin having the maximum
contra-harmonic mean could potentially be selected for diving. Similarly, when shifting priority to the small items, harmonic means of
the bins could potentially be considered as the diving criterion.

An alternative diving criterion could further be obtained by considering different Lehmer
means simultaneously. The criterion
\begin{equation}
\label{integral}
L_s = \int_{0}^{2} p L_p(W)  dp \ ,
\end{equation}
assumes different Lehmer means ranging from harmonic to contra-harmonic
means where the impacts of the smaller items are relaxed by applying the
multiplier $p$. The following example illustrates an application of the geometrical
diving criteria.

\begin{figure}[!t]
        \centering
        \includegraphics[width=\linewidth, height=6cm]{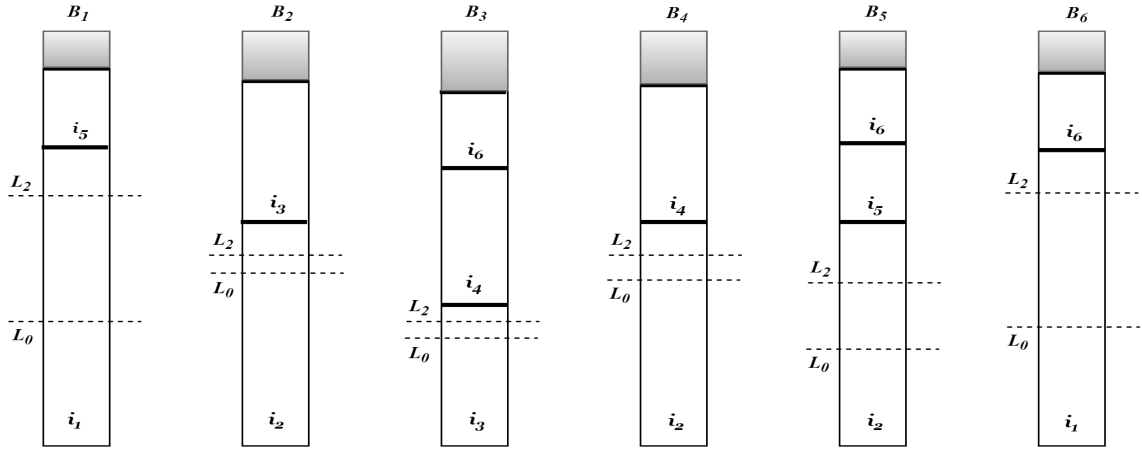}		        
        \caption{Composition of items into bins of Example \ref{Example}.}
		\label{fig:Example}  
\end{figure}

\begin{example}
\label{Example}
Let $n = 6$, $c = 100$, $I = \{1, 2, 3, 4, 5, 6\}$, $W_I=\{72, 54, 34, 33, 19, 18\}$
and $D_I=\{1, 1, 1, 1, 1, 1\}$. Solving the master problem at the root node,
yields the following set of bins
$$B_\mathrm{rel} = \left\{\{i_1 , i_5\}, \{i_2 , i_3\}, \{i_3 , i_4 , i_6 \}, \{i_2 , i_4 \}, \{i_2 , i_5 , i_6 \}, \{i_1 , i_6 \} \right\} \ , $$
where the relaxed values bins are obtained as being $0.8, 0.4, 0.6, 0.4, 0.2$ and $0.2$, respectively.
By setting the diving criterion to be the highest-value, diving at the root
node occurs on $B_1 = \{i_1 , i_5\}$. However, by setting the criterion to be $L_0$ ,
$B_2 = \{i_2 , i_3\}$ gets selected whose small item is larger than the small items
contained in the other bins. On the other hand, the criterion $L_2$ selects
$B_6 = \{i_1 , i_6 \}$ where the larger item $i_1$ is prioritized. Furthermore, by evaluating
the integral (\ref{integral}) for each bin, $B_1$ gets selected as the candidate for diving.

For this simple example, employing any of the geometrical diving criteria leads to the primal heuristic
solution having the optimal number of bins. Composition of the items into
bins in $B_\mathrm{rel}$ and the corresponding harmonic and contra-harmonic means
for each bin are depicted in Fig. \ref{fig:Example}.

\end{example}

\section{Batch Diving}
\label{BatchDivingSec}
Practically speaking, the diving methods described in Sect. \ref{DivingSec} could not be
implemented efficiently since the computational time for re-optimization of the child nodes remain close to the computational time of their parent nodes.
For this reason, a batch diving procedure is proposed where diving occurs
on a batch of bins concurrently.

Let us assume that $\bar{n}$ items are present in a certain node ($\bar{n}=n$ at the root
node) and let $\bar{d}_i$ denote the corresponding residual demand of each item, where $i = 1, \dots , \bar{n}$. The model for batch diving is introduced as
\begin{subequations} 
\label{BD_obj}
\begin{alignat}{3}
\label{BD_1}
\min & \ \sum_{k=1}^{K^\star} \gamma_k \mu_k  \\
\label{BD_2}
\mathrm{s.t.} & \ \sum_{k=1}^{K^\star} x_{ik} \mu_k= \bar{d}_i && \quad i=1,\dots,\bar{n} \\
 & \ \mu_k \in \{0,1\} &&\quad k=1,\dots,K^\star \ ,
\end{alignat}
\end{subequations}
where $K^\star$ is the number of bins obtained when master problem is solved at
the considered node. Solving (\ref{BD_obj}) returns a set of high quality bins,
with respect to the diving criterion. In (\ref{BD_1}), $\gamma_k$ is the value of the diving
criterion for bin $k$ and $\mu_k$ is the load of bin $k$, where $\mu_k = 1$ if the $k$-th bin is present in the
set of bins candidate for batch diving, $0$ otherwise. Also, the set of constraints
(\ref{BD_2}) ensure that diving on a batch of some selected bins would not lead to
an infeasible node.

After batch diving is performed at each node, dimensions of the residual
problem decrease significantly, which in turn would speed up the re-optimization
process. It is also possible to perform the batch diving on smaller subsets of
bins at each node in order to increases the chances of constructing the optimal solution to BPP.

\section{Polluted Pricing Sub-problem}
The branching rule in the \texttt{GEOM-BP} algorithm is the very same criterion
used in diving, and more attention is put on forbidden bins by exploring the 
tree in a binary manner. As a consequence, a single bin is added either to
the partial upper bound solution or to the list of the forbidden bins. In the
binary, compared to the integer branching strategy, more nodes containing
forbidden bins are explored. It is worth mentioning that the paths to
optimal solutions of some of the most difficult instances of BPP are likely to be the ones that visit the polluted nodes.

As mentioned earlier, the main concern in developing branch-and-price
algorithms pertains to overcoming the issue of the presence of forbidden
bins polluting the search space of the pricing sub-problem. Our approach
when dealing with the polluted pricing sub-problem of column generation
is straightforward. In practice, whenever one of the forbidden bins is met,
a two-dimensional knapsack problem (2D-KSP) is solved in order to derive
the next feasible solution to the pricing sub-problem.
Assume that a forbidden bin is generated at the $k$-th iteration of column
generation with an objective value $\nu_k$ where
\begin{equation}
\nu_k = \sum\limits_{i=1}^{n} \pi_i x_i - 1 \ ,
\end{equation}
then the next feasible solution to the pricing sub-problem could be obtained by solving 
\begin{subequations} 
\label{2D_KS}
\begin{alignat}{3}
\label{2D_KS_1}
\nu_{\mathrm{2D-KSP}} (\boldsymbol{\pi}) = \max & \ \sum_{i=1}^{n} \pi_i x_i - 1  \\
\label{2D_KS_2}
\mathrm{s.t.} & \ \sum_{i=1}^{n} \pi_i x_i \leq c &&  \\
\label{2D_KS_3}
 & \ \sum_{i=1}^{n} \pi_i x_i - 1 \leq \nu_k - \delta &&  \\
 & \ x_i \in \{0,1,\dots,\bar{d}_i\} &&\quad i=1,\dots,\bar{n} \ ,
\end{alignat}
\end{subequations}
where $\delta$ denotes the value of a decrement. In this model, (\ref{2D_KS_3}) is referred to as 
the decrement constraint, and it ensures that the generated bin will have an
objective value lower than $\nu_k$. In other words, (\ref{2D_KS}) guarantees that
an arbitrary generated bin contains a composition of items different from
that of the forbidden one. The process is repeated in case the generated bin
after solving (\ref{2D_KS}) appears to be one of the other forbidden bins again.

An important aspect of (\ref{2D_KS}) is the value assigned to the decrement.
For large values of $\delta$, there could emerge possibilities to skip some high
quality solutions of the sub-problem. More importantly, the bin generated by this
model could cause termination of column generation. In such a case, since
solving (\ref{2D_KS}) does not guarantee the achievement of the second best
solution of the sub-problem, column generation might terminate prematurely.

\begin{figure}[!t]
        \centering
        \includegraphics[width=\linewidth, height=10cm]{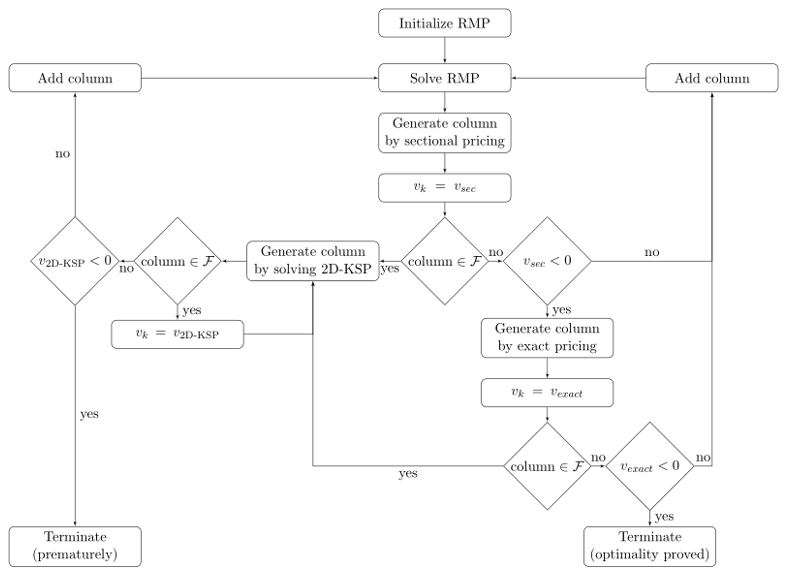}
        
        \caption{Schematic of steps for solving master problems in \texttt{GEOM-BP}.}
        \label{fig:Flowchart}

\end{figure}

It is noted that even when considering a minimal decrement, achievement
of the second best solution is not guaranteed. For these nodes, pruning by
bound is not allowed and branching is resumed. Usually, by considering
sufficiently small values for the decrement $\delta$, optimality of the column generation gets successfully proved for the majority of the nodes. Figure \ref{fig:Flowchart} depicts an schematic of the steps described above where $\mathcal{F}$ denotes the set of forbidden bins.

\section{Computational Results}
In this section, we present computational results for running the \texttt{GEOM-BP} algorithm on the benchmark instances. Armadillo C\texttt{++} \cite{sanderson2016armadillo} library was
used as the main hub for programming, and the C codes for \texttt{COMBO} \cite{martello1999dynamic} and \texttt{BOUKNAP} \cite{pisinger2000minimal} were compiled in GCC 7.2.0 to solve the binary and bounded 1D-KSP, respectively. Furthermore, the 2D-KSP were solved calling the Fortran sub-routine provided by Martello et al. \cite{martello2003exact}, and the multi-dimensional knapsack problems arising in the batch diving models were solved employing CPLEX 12.7.0 . All experiments
were carried out on an Intel core i-7 6700 HQ 2.6 GHz with 16 GB RAM,
but only a single core was allowed when executing the algorithms to have a
reasonable comparison with other methods.

The classes of benchmark instances considered in this work include:
\begin{itemize}
\item A set of 80 instances `Falkenauer U' and a set of 80 triplets `Falkenauer
T' presented by Falkenauer in \cite{falkenauer1996hybrid}.
\item Sets of 720 instances `Scholl 1', 480 instances `Scholl 2' and 10 instances
`Scholl 3' presented by Scholl et al. in \cite{scholl1997bison}.
\item A set of 17 instances `W{\"a}scher' presented by W{\"a}scher and Gau in \cite{wascher1996heuristics}.
\item Sets of 100 instances `Schwerin 1' and 100 instances `Schwerin 2' presented by Schwerin and W{\"a}scher in \cite{schwerin1997bin}.
\item A set of 28 instances `Hard28' presented by Schoenfield in \cite{schoenfield2002fast}.
\item A set of 3840 instances `Delorme R' presented by Delorme et al. in \cite{delorme2016bin}.
\end{itemize}
However, the trivial instances from the literature whose computed combinatorial lower bounds coincide with the number of bins resulting from the Best Fit Decreasing algorithm were excluded from experiments. The combinatorial lower bound for each instance could be computed by a
method presented by Martello and Toth \cite{martello1990lower}. 

Table \ref{CompTable} compares the results obtained using \texttt{GEOM-BP} to the exact algorithms from the literature, namely \texttt{MTP}\cite{martello1990lower}, \texttt{BISON} \cite{scholl1997bison}, \texttt{ONECUT} \cite{dyckhoff1981new}, \texttt{ARCFLOW} \cite{de1999exact}, \texttt{VPSOLVER} \cite{brandao2016bin}, \texttt{VANCE} \cite{vance1994solving} and \texttt{BELOV} \cite{belov2006branch}. We adapt the computational results reported by Delorme et al. \cite{delorme2016bin} as the bechmark. It is observed that Geom-BP algorithm, using any of the introduced geometrical criteria $L_0$, $L_2$ or $L_s$ outperforms the previously developed exact algorithms in the tested number of instances solved in less than one minute. 

It is noteworthy to mention that \texttt{GEOM-BP} shows a significant
improvement over BELOV algorithm when solving the instances of class
`Falkenauer T'. This is
mainly due to the fact that these instances hold Round-Up Property (RUP)
where the objective value of the master problem solved at the root node is
equal to the optimal number of bins, and the use of our proposed heuristic
together with the batch diving procedure prove useful in constructing the
optimal bins as fast as possible. 

For the more challenging instances
of class `Hard28', \texttt{BELOV} algorithm remains superior both in the number
of instances solved in less than one minute and the average computational
time. This class contains $5$ instances with Non-Integer Round-Up Property (Non-IRUP) where the gap between the integer solution and the LP relaxation at the root node is greater than or equal to $1$. For these Non-IRUP instances, the gaps are not closed and \texttt{GEOM-BP} proves the optimality of the solutions by exploring and exhausting all the possible nodes. It is pointed out that an appropriate choice for the decrement value $\delta$ in (\ref{2D_KS}) leads to termination of branch-and-price algorithm in less than one minute. In this work, we use an adaptive scheme to determine the decrement value whose initial point is considered to be the relatively small value $\delta_0 = 10^{-5}$. For 2 instances from the `Hard28' class, namely `BPP60' and `BPP181', optimality of the solutions were proved only by increasing the time limit to $10$ minutes. However, \texttt{GEOM-BP} could not close the gap for the instance `BPP40' in less than $10$ minutes, and this instance remained unsolved. 

For the `Scholl 3' instances, the bottleneck for solving the master problem at the root node pertains to the large dimenstions of the sub-problems, preventing the competitive algorithms like \texttt{VPSOLVER} and \texttt{BELOV} from solving the instances efficiently. The sectional pricing scheme proposed in this work, leads \texttt{GEOM-BP} to solve these instances in shorter amounts of time by reducing the dimensions of the pricing sub-problems.

Further details for \texttt{GEOM-BP} are presented in Table (\ref{DetailTable}).
For each class of instances, the average total number of nodes explored
($n_{\mathrm{total}}^{\mathrm{node}}$) and the average total number of polluted nodes explored ($n_{\mathrm{poll}}^{\mathrm{node}}$)
are reported. It is seen that for majority of the instances, the optimal solution is obtained without exploring the polluted nodes. However, for more
difficult classes of `Falkenauer T', `W{\"a}scher' and `Hard28', optimal solutions
are achieved only when the polluted nodes are explored. Furthermore, $n_{\mathrm{col}}^{\mathrm{root}}$ denotes the average number of the columns generated at the root node and $n_{\mathrm{exact}}^{\mathrm{root}}$ represents average number of times the exact pricing scheme is called to prove the optimality of column generation. Even though this number increases when solving some
instances, especially the ones from the class of `W{\"a}scher', still the overall
computational time when using the implicit sectional pricing scheme is less
than that of the only exact pricing scheme.


\begin{landscape}
\begin{table}[!t]
\centering
\caption{Number of tested instances (average time in seconds) solved in less than one minute}
\label{CompTable}
\begin{tabular}{@{}lrrrrrrrrr@{}}
\toprule
Class         & Inst. & MTP         & BISON       & ONECUT     & ARCFLOW    & VPSOLVER   & VANCE       & BELOV      & GEOM-BP    \\ \midrule
Falkenauer U  & 74    & 22 (42.8)   & 44 (24.5)   & 74 (0.2)   & 74 (0.2)   & 74 (0.1)   & 53 (24.1)   & 74 (0.0)   & \textbf{74} (0.0)   \\
Falkenauer T  & 80    & 6 (55.5)    & 42 (30.6)   & 80 (8.7)   & 80 (3.5)   & 80 (0.4)   & 76 (14.8)   & 57 (24.7)  & \textbf{80} (0.2)   \\
Scholl 1      & 323   & 242 (15.1)  & 288 (7.0)   & 323 (0.1)  & 323 (0.1)  & 323 (0.1)  & 323 (3.6)   & 323 (0.0)  & \textbf{323} (0.0)  \\
Scholl 2      & 244   & 130 (28.2)  & 233 (3.0)   & 118 (38.7) & 202 (18.9) & 208 (14.0) & 204 (18.6)  & 244 (0.3)  & \textbf{244} (0.1)  \\
Scholl 3      & 10    & 0 (60.0)    & 3 (42.0)    & 0 (63.9)   & 0 (61.1)   & 10 (6.3)   & 10 (1.9)    & 10 (14.1)  & \textbf{10} (5.3)   \\
W{\"a}scher & 17    & 0 (60.0)    & 10 (24.7)   & 0 (60.7)   & 0 (60.5)   & 6 (49.4)   & 6 (52.0)    & 17 (0.1)   & \textbf{17} (0.2)   \\
Schwerin 1    & 100   & 15 (51.1)   & 100 (0.0)   & 100 (13.1) & 100 (1.5)  & 100 (0.3)  & 100 (0.3)   & 100 (1.0)  & \textbf{100} (0.0)  \\
Schwerin 2    & 100   & 4 (57.6)    & 63 (22.2)   & 100 (11.7) & 100 (1.5)  & 100 (0.3)  & 100 (0.3)   & 100 (1.4)  & \textbf{100} (0.0)  \\
Hard28        & 28    & 0 (60.0)    & 0 (60.0)    & 6 (54.6)   & 16 (40.6)  & 27 (14.2)  & 11 (48.9)   & 28 (7.3)   & 25 (11.2)  \\
Delorme R     & 2901  & 1252 (34.7) & 1839 (22.6) & 2789 (5.0) & 2840 (3.3) & 2891 (1.4) & 2152 (20.3) & 2901 (0.2) & 2898 (0.7) \\ \midrule
Overall       & 3877  & 1671 (34.6) & 2652 (20.0) & 3590 (7.8) & 3735 (4.5) & 3819 (2.3) & 3035 (18.0) & 3854 (0.8) & \textbf{3871} (0.6) \\ \bottomrule
\end{tabular}
\end{table}

\begin{table}[!t]
\centering
\caption{Computational details for \texttt{GEOM-BP}}
\label{DetailTable}
\begin{tabular}{lrrrrr}
\hline
Class         & Inst. & $n_{\mathrm{col}}^{\mathrm{root}}$    & $n_{\mathrm{exact}}^{\mathrm{root}}$ & $n_{\mathrm{total}}^{\mathrm{node}}$ & $n_{\mathrm{poll}}^{\mathrm{node}}$ \\ \hline
Falkenauer U  & 74    & 353.9  & 1.5      & 5.7      & 0.0     \\
Falkenauer T  & 80    & 899.8  & 1.4      & 14.7     & 0.2     \\
Scholl 1      & 323   & 254.6  & 1.6      & 4.4      & 0.0     \\
Scholl 2      & 244   & 3041.6 & 1.2      & 4.8      & 0.0     \\
Scholl 3      & 10    & 1587.1 & 1.0      & 4.4      & 0.0     \\
W{\"a}scher & 17    & 740.4  & 6.8      & 9.2      & 3.4     \\
Schwerin 1    & 100   & 198.4  & 1.1      & 3.6      & 0.0     \\
Schwerin 2    & 100   & 269.8  & 1.1      & 5.2      & 0.0     \\
Hard28        & 28    & 7419.8 & 1.4      & 82.0     & 15.3    \\
Delorme R     & 2901  & 1258.1 & 1.2      & 6.3      & 0.0     \\ \hline
Overall       & 3877  & 1352.7 & 1.4      & 6.3      & 0.1     \\ \hline
\end{tabular}
\end{table}

\end{landscape}

\bibliographystyle{splncs}      
\bibliography{GeomBP}            

\end{document}